\documentclass[a4paper, 11pt, oneside, final]{scrartcl}

\usepackage{PREAMBLE}
\usepackage{COMMAND}
\usepackage{HYPHENATION}

\begin{document}
\thispagestyle{empty}
\begin{center}
    \rule{\linewidth}{1pt}\\[0.4cm]
    {\sffamily \bfseries \large Enforcing uniqueness in one-dimensional phase retrieval by
        additional signal information in time domain}\\[10pt]
    {\sffamily\footnotesize Robert Beinert and Gerlind Plonka}\\[3pt]
    {\sffamily\footnotesize Institut für Numerische und Angewandte Mathematik}\\[-3pt]
    {\sffamily\footnotesize Georg-August-Universität Göttingen}\\
    \rule{\linewidth}{1pt}
\end{center}

\vspace*{10pt}

\begin{center}
    \sffamily\bfseries Abstract
\end{center}

{\small
    \noindent
    Considering the ambiguousness of the discrete-time phase retrieval problem to recover
    a signal from its Fourier intensities, one can ask the question: what additional
    information about the unknown signal do we need to select the correct solution within
    the large solution set?  Based on a characterization of the occurring ambiguities, we
    investigate different a priori conditions in order to reduce the number of ambiguities
    or even to receive a unique solution.  Particularly, if we have access to additional
    magnitudes of the unknown signal in the time domain, we can show that almost all
    signals with finite support can be uniquely recovered.  Moreover, we prove that an
    analogous result can be obtained by exploiting additional phase information.

    \medskip


    \smallskip
    
}

\section{Introduction}
\label{sec:introduction}

The phase retrieval problem consists in recovering a complex-valued signal from the
modulus of its Fourier transform.  In other words, the phase of the signal in the
frequency domain is lost.  Recovery problems of this kind have many applications in
physics and engineering as for example in crystallography \cite{Mil90,Hau91,KH91},
astronomy \cite{DF87}, and laser optics \cite{SST04,SSD+06}.  Finding an analytic or a
numerical solution is generally challenging due to the well-known ambiguousness of the
problem.  In order to determine a meaningful solution, one hence requires further
appropriate information about the unknown signal.

In this paper, we consider the phase retrieval problem in a discrete setting and restrict
ourselves to the recovery of a complex-valued discrete-time signal
$x \coloneqq (x[n])_{n\in\Z}$ with finite support form its Fourier intensities
$\absn{\fourier x}$.  All occurring ambiguities of this problem can be explicitly
constructed via the zeros of the autocorrelation polynomial, see \cite{BS79,OS89,BP15}.

Depending on the application and the exact problem setting, different a priori conditions
have been employed in the literature to extract special solutions of the phase retrieval
problem.  For real-valued signals, the interference with a known reference signal
\cite{KH90a,KH90} can be exploited to reduce the complete solution set to at most two
different signals.  Further, the interference with an unknown reference signal has been
considered in \cite{KH93}.  These approaches can also be used to recover complex-valued
signals, see \cite{RDN13,RSA+11,BP15}.  Moreover, comparable results have been achieved by
special reference signals that are strongly related to the unknown signal itself.  For
example, one can use a modulated version of the original signal, see \cite{CESV13,ABFM14}.

Other approaches are based on additional information in the frequency domain.  In
\cite{HHLO83}, beside the modulus $\absn{\fourier x}$, the information whether the phase of Fourier
transform $\fourier x (\omega)$ is contained in
$[-\nicefrac{\pi}{2}, \nicefrac{\pi}{2}]$ or in
$[-\pi, - \nicefrac{\pi}{2}) \cup (\nicefrac{\pi}{2}, \pi)$ has been studied.  One
may also replace the Fourier transform by the so-called short-time Fourier transform
\cite{NQL83a,NQL83}, where the original signal is overlapped with a small analysis window
at different positions.

If the unknown discrete-time signal has a fixed support of the form $\{0,\dots,M-1\}$ and can thus be
identified with $M$-dimensional vectors, then the Fourier intensities
$\absn{\fourier x(\omega_k)}$ at different points $\omega_k \in [-\pi,\pi)$ can be
written as the intensity measurement $\absn{\iProd{x}{v_{k}}}$ with
$v_{k} \coloneqq (\e^{\I \omega_k m})_{m=0}^{M-1}$.  In a more general setting using a
frame approach, the question arises how the vectors $v_k$ have to be constructed, and how
many vectors $v_k$ are needed to ensure a unique recovery of $x$ only from the intensities
$\absn{\iProd{x}{v_k}}$, see for instance \cite{BCE06, BBCE09, ABFM14,BCM14,BH15} and
references therein.

\smallskip

In some applications, like wave front sensing and laser optics \cite{SST04}, besides the
Fourier intensity, the modulus of the unknown signal itself is known.  For this specific
one-dimensional phase retrieval problems, a multi-level Gauss-Newton method has been
presented in \cite{SSD+06,LT08,LT09} to determine a numerical solution.  While this method
worked well for the considered problems, its stability seems to depend on the given data
sets \cite{Bei13a}.  Further, for some rare cases, the algorithm converges to an
approximate solution which is completely different from the original signal.

The occurring numerical problems in the above mentioned multi-level Gauss-Newton method
have been the reason for our extensive studies of ambiguities of the complex discrete
phase retrieval problem within this paper.  Particularly, we will consider the question
whether the additional knowledge of $\absn{x} \coloneqq (\absn{x[n]})_{n\in \Z}$ can
indeed ensure a unique recovery of the unknown signal $x$. When we understand this
uniqueness problem completely, then we can decide whether the behavior of the algorithm is
due to ambiguities or just to ill-posedness of the problem.  Besides modulus constraints
in time domain, we will also consider the case that information about the phase
$\arg x \coloneqq (\arg x[n])_{n\in \Z}$ in time domain is available, and investigate how
far this additional phase information can reduce the solution set.

\smallskip

The paper is organized as follows. In section~\ref{sec:phase-retr-probl}, we briefly recall
the characterization of the solution set in \cite[Theorem~2.4]{BP15}, which is based on an
appropriate factorization of the autocorrelation function.  Distinguishing the trivial
ambiguities, like rotations, shifts, and conjugations and reflections, from the
non-trivial, we will observe that the discrete-time phase retrieval problem possesses only
a finite set of relevant solutions.

In section~\ref{sec:using-add-magn}, we investigate the recovery of a complex-valued signal
$x$ with fixed support $\{0, \dots, N-1\}$ from its Fourier intensity $\absn{\fourier x}$
where the modulus $\absn{x[n]}$ of at least one signal component in time domain is also
given.  For that purpose, we generalize the findings in \cite{XYC87,BP15} that almost each
signal $x$ is uniquely determined by its Fourier transform $\absn{\fourier x}$ and its end
point $x[N-1]$.  We will show that already one given modulus $\absn{x[n]}$ of the unknown
signal can enforce uniqueness (up to trivial ambiguities) for almost every signal.
Unfortunately, the uniqueness cannot be achieved for every signal even if we know the
complete modulus $\absn{x}$. We will construct examples where uniqueness is not obtained.

Finally, in section~\ref{sec:add-phase-point}, we consider the phase retrieval problem with
additional phase information $\arg x[n]$ at appropriate points in the time domain and
prove that already two given phases can avoid the ambiguousness for almost all signals.
However, similarly to the case of given moduli, even the complete phase information
$\arg x \coloneqq (\arg x[n])_{n\in\Z}$ may be not sufficient to ensure the uniqueness for
special cases.

\section{The phase retrieval problem}
\label{sec:phase-retr-probl}

We consider the \emph{one-dimensional discrete-time phase retrieval problem} where we wish
to recover the complex-valued discrete-time signal $x \coloneqq (x[n])_{n\in\Z}$ with
finite support from its Fourier intensity
\begin{equation*}
    \absn{\fourier x \mleft( \omega \mright) }
    \coloneqq \abs{\Fourier \mleft[ x \mright] \mleft( \omega \mright)}
    \coloneqq \absB{\sum_{n\in\Z} x\mleft[n\mright] \, \e^{-\I \omega n}}
    \qquad (\omega \in \R).
\end{equation*}
Unfortunately, this problem is complicated because of the well-known ambiguousness.  For
example, we can easily construct further solutions by rotating, shifting, or reflecting
and conjugating the original signal $x$, see \cite{BP15}.

\begin{Proposition}\label{prop:trivial-amb}
    Let $x$ be a complex-valued signal with finite support. Then
    \begin{enumerate}[\upshape(i)]
    \item\label{item:1} the rotated signal $( \e^{\I \alpha} \, x[ n] )_{n\in\Z}$ for
        $\alpha \in \R$
    \item\label{item:2} the time shifted signal $( x[n-n_0] )_{n\in\Z}$ for $n_0 \in \Z$
    \item\label{item:3} the conjugated and reflected signal $(\overline{ x[-n]})_{n\in\Z}$
    \end{enumerate}
    have the same Fourier intensity $\absn{\fourier x}$.%
\end{Proposition}

\begin{Proof}
    This observation can be simply verified by determining the Fourier transforms of the
    considered signals.  Obviously, we have
    \begin{enumerate}[(i)]
    \item $\Fourier[(\e^{\I \alpha} \, x[n])_{n\in \Z}] = \e^{\I \alpha} \, \fourier x;$
    \item $\Fourier[(x[n - n_0])_{n\in\Z}] = \e^{-\I \omega n_0} \, \fourier x;$
    \item $\Fourier[(\overline{x[ -n ]})_{n\in\Z}] = \overline{\fourier x}.$ \qed
    \end{enumerate}
\end{Proof}

Although the rotation, shift, and conjugation and reflection of a solution always results
in a further solution, these signals are very closely related to the original signal.
Therefore, we call these three kinds of ambiguities \emph{trivial} as introduced in
\cite{Wan13}. However, besides the trivial ambiguities, the discrete-time phase retrieval
problem usually possesses an extensive amount of further non-trivial solutions, see
\cite{BS79,BP15}.

We briefly recall an explicit characterization of the complete solution set, which is
based on the observations in \cite{BP15}.  Let the \emph{autocorrelation signal} $a$ of a
signal $x$ be defined by
\begin{equation*}
    a \mleft[ n \mright] \coloneqq
    \sum_{k\in\Z} \overline{x \mleft[ k \mright]} \, x \mleft[ k+n \mright]
    \qquad
    (n \in \Z).
\end{equation*}
Further, let the \emph{autocorrelation function} be the Fourier transform $\fourier a$ of
the autocorrelation signal $a$.  If $N$ denotes the support length of the signal $x$, then
the corresponding autocorrelation signal $a$ possesses the support
$\{-N+1, \dots, N-1 \}$.  Moreover, the autocorrelation signal $a$ is always conjugate
symmetric, i.e., $a[-n] = \overline{a[n]}$ for $n \in \Z$. Thus the autocorrelation
function $\fourier a$ is a real-valued trigonometric polynomial of degree $N-1$ of the
from
\begin{equation*}
    \fourier a \mleft( \omega \mright)
    = \sum_{-N+1}^{N-1} a \mleft[ n \mright] \, \e^{-\I \omega n}
    \qquad(\omega \in \R).
    \addmathskip
\end{equation*}
Moreover, we have the relation
\begin{equation*}
    \fourier a \mleft( \omega \mright)
    = \sum_{n\in \Z} \sum_{k \in \Z} x\mleft[k+n\mright] \, \overline{x\mleft[k\mright]} \, \e^{-\I\omega n}
    = \sum_{n\in \Z} \sum_{k\in\Z} x \mleft[n\mright] \, \overline{x \mleft[k\mright]} \,
    \e^{-\I \omega (n-k)}
    = \absn{\fourier x \mleft( \omega \mright)}^2\!.
\end{equation*}
Thus, knowing the Fourier intensity $\absn{\fourier x}$ is equivalent to knowing the
autocorrelation function.  Further, the real-valued, non-negative trigonometric polynomial
$\fourier a$ of degree $N-1$ is already determined by $2N-1$ samples from
$[-\pi,\pi)$.  Consequently, our initial phase retrieval problem is equivalent to the
completely discrete formulation: Recover the complex-valued signal
$x \coloneqq (x[n])_{n\in\Z}$ with support of length $N$ from the $2N-1$ values
\begin{equation*}
    \abs{\fourier x \mleft( \tfrac{2\pi k}{N} \mright)}
    \qquad
    (k=-N,\dots,N-1).
    \addmathskip
\end{equation*}

We summarize the considerations in \cite{BP15} to characterize all solutions of
$\absn{\fourier x}^2 = \fourier a$ explicitly.  Let $P(z)$, given by
$ P( \e^{-\I \omega} ) = \e^{-\I \omega (N-1)} \, \fourier a ( \omega )$, be the
\emph{associated polynomial} to $ \fourier a$, i.e.
\begin{equation*}
    P \mleft( z \mright)
    \coloneqq \sum_{n=0}^{\crampedclap{2N-2}} a \mleft[n-N+1 \mright] \, z^n.
    \addmathskip
\end{equation*}
Here the conjugate symmetry $a[-n] = \overline{a[n]}$ of the coefficients implies that the
zeros of the associated polynomial $P$ appear in reflected pairs
$(\gamma_j^{}, \overline \gamma_j^{\,-1})$ with respect to the unit circle, and that the
zeros on the unit circle have even multiplicity.  Since $P$ is a polynomial of degree
$2N-2$, there are precisely $N-1$ such zero pairs, yielding the factorization
\begin{equation}
    \label{eq:ass-poly-fact}
    P \mleft( z \mright)
    = a \mleft[ N-1 \mright] \, \prod_{j=1}^{N-1} \left( z - \gamma_j^{\,} \right) \left(
        z - \overline \gamma_j^{\,-1} \right) \!.
\end{equation}
Now, we can employ \eqref{eq:ass-poly-fact} to factorize $\fourier a$,  
\begin{equation*}
    \fourier a \mleft( \omega \mright) = \abs{\fourier a \mleft( \omega \mright)}
    = \abs{ P \mleft( \e^{-\I \omega} \mright)}  
    = \abs{a \mleft[ N-1 \mright]} \prod_{j=1}^{N-1} \abs{\e^{-\I\omega} -
        \gamma_j^{\,}} \abs{\e^{-\I\omega} - \overline \gamma_j^{\,-1}} \!.
\end{equation*}
We rewrite the second linear factor by using the identity
\begin{equation*}
    \absb{\e^{-\I\omega} - \overline \gamma_j^{\,-1}}
    = \absn{\overline \gamma_j}^{-1} \, \absb{\overline \gamma_j - \e^{\I\omega}}
    = \absn{\gamma_j^{\,}}^{-1} \, \absb{\e^{-\I \omega} - \gamma_j^{\,}} 
\end{equation*}
and have the representation
\begin{equation}
    \label{eq:autocor-func-fact}
    \fourier a \mleft( \omega \mright)
    = \abs{ a \mleft[ N-1 \mright]} \prod_{j=1}^{N-1} \absn{\gamma_j}^{-1}
    \cdot \absB{ \prod_{j=1}^{N-1} \left( \e^{-\I \omega} - \gamma_j \right)}^2.
\end{equation}
Since one can similarly rewrite the first linear factor instead of the second, we obtain
the following characterization of the solution set of the discrete-time phase retrieval
problem, see \cite[Theorem~2.4]{BP15}.

\begin{Theorem}
    \label{the:repr-sol-time-dom}%
    Let $\fourier a$ be a non-negative trigonometric polynomial of degree $N-1$.  Then,
    each solution $x$ of the discrete phase retrieval problem
    $\absn{\fourier x}^{2} = \fourier a$ with finite support has a Fourier representation
    of the form
    \begin{equation}
        \label{eq:char-amb}
        \fourier x \mleft( \omega \mright)
        = \e^{\I \alpha} \, \e^{-\I\omega n_0}
        \sqrt{\abs{a \mleft[ N-1 \mright]} \prod_{j=1}^{N-1} \abs{\beta_j}^{-1}} \cdot
        \prod_{j=1}^{N-1} \left( \e^{-\I \omega} - \beta_j \right) \!,
        \addmathskip
    \end{equation}
    where $\alpha$ is a real number, $n_0$ is an integer, and where for each $j$ the value
    $\beta_j$ is chosen from the zero pair $(\gamma_j^{\,}, \overline \gamma_j^{\,-1})$ of
    the associated polynomial to $\fourier a$.%
\end{Theorem}

Apart from the trivial rotation and shift ambiguity caused by the unimodular factor
$\e^{\I\alpha}$ and the modulation $\e^{-\I \omega n_0}$ in Theorem~\ref{the:repr-sol-time-dom},
each solution $y$ is characterized by the values
$\beta_j \in (\gamma_j^{\,}, \overline \gamma_j^{\,-1})$.  Hence, we can uniquely identify
each solution $x$ up to rotations and shifts with the set
$B \coloneqq \{ \beta_1, \dots, \beta_{N-1} \}$, which we call the \emph{corresponding
    zero set} of the solution $x$ in the following. Obviously, we can construct up to
$2^{N-1}$ different zero sets; the corresponding signals, however, do not have to be
non-trivially different, see \cite[Corollary~2.6]{BP15}.

\begin{Lemma}
    \label{lem:zeros-refl-sig}
    Let $x$ be a discrete-time signal of the form \eqref{eq:char-amb} with corresponding
    zero set $B \coloneqq \{ \beta_1, \dots, \beta_{N-1} \}$.  Then, the conjugated and
    reflected signal $\overline{x [- \cdot ]}$ corresponds to the zero set
    \begin{equation*}
        \bigl\{ \overline \beta_1^{\,-1}, \dots, \overline \beta_{N-1}^{\,-1} \bigr\} .
    \end{equation*}
\end{Lemma}

Thus we can conclude that the phase retrieval problem to recover $x$ can possesses up to
$2^{N-2}$ non-trivially different solutions, \cite[Corollary~2.6]{BP15}.  The number
$2^{N-2}$ is here only an upper bound for the occurring non-trivial solutions since the
actual number for a specific phase retrieval problem strongly depends on the zero pairs
\raisebox{0pt}[0pt][0pt]{$(\gamma_j^{\,}, \overline \gamma_j^{\,-1})$} of the associated
polynomial.  For example, if all zero pairs
\raisebox{0pt}[0pt][0pt]{$(\gamma_j^{\,}, \overline \gamma_j^{\,-1})$} lie on the unit
circle, then every $\beta_j$ in Theorem~\ref{the:repr-sol-time-dom} is uniquely determined and
the corresponding phase retrieval problem is uniquely solvable.

\section{Using additional magnitudes of the unknown signal}
\label{sec:using-add-magn}

To determine a unique solution within the solution set characterized by
Theorem~\ref{the:repr-sol-time-dom}, we need additional information about the unknown signal.
In \cite{BP15}, we have already shown that almost every signal $x$ with support
$\{0, \dots, N-1 \}$ is uniquely determined by its Fourier intensity $\absn{\fourier x}$
and the right end point $x[N-1]$.  A similar observation has been done for real-valued
signals by Xu et al.\ in \cite{XYC87}.  In this section, we will show that comparable
results can be achieved by exploiting one absolute value $\absn{x[n]}$ or even several
absolute values $\absn{x[n]}$ of the unknown signal $x$.

\subsection{The modulus of an arbitrary signal value}
\label{sec:modulus-arb-point}%

We consider the phase retrieval problem to recover the signal $x$ with support
$\{0, \dots, N-1\}$ from its Fourier intensity $\absn{\fourier x}$ and the modulus of one
signal value $\absn{x[N-1-\ell]}$ for some $\ell$ between $0$ and $N-1$.  We will show
that this phase retrieval problem is almost surely uniquely solvable up to rotations.

Let $B \coloneqq \{ \beta_1, \dots, \beta_{N-1} \}$ be the corresponding zero set of a solution
signal $x$ as given in (\ref{eq:char-amb}) and the subsequent comments. Further, let
$\Lambda$ be a subset of $B$.  We introduce the \emph{modified zero set}
\begin{equation}\label{B+}
    B^{(\Lambda)}
    \coloneqq
    \left\{\beta_1^{(\Lambda)}, \dots, \beta_{N-1}^{(\Lambda)} \right\},
    \addmathskip
\end{equation}
where the single elements are given by
\begin{equation*}
    \beta_j^{(\Lambda)} \coloneqq
    \begin{cases}
        \overline \beta_j^{\,-1} & \beta_j \in \Lambda, \\
        \beta_j & \text{else}.
    \end{cases}
\end{equation*}
Recall that by Theorem \ref{the:repr-sol-time-dom} each further non-trivial solution $y$
of the discrete-time phase retrieval problem $\absn{\fourier x}^2 = \fourier a$
corresponds to such a modified zero set.

Let the $(N-1)$-variate \emph{elementary symmetric polynomial} $S_n: {\C}^{N-1} \to {\C}$
of degree $n$ in the variables $\beta_1$, \dots, $\beta_{N-1}$ be given by
\begin{equation}
    \label{eq:ele-sys-fun}
    S_n \mleft( \beta_1, \dots, \beta_{N-1} \mright)
    \coloneqq
    \smashoperator{\sum_{\quad1 \le k_1 < \dots < k_n \le N-1}} \beta_{k_1} \cdots \beta_{k_n}
\end{equation}
for $n$ from $1$ to $N-1$.  Since the polynomials $S_n$ are independent from the order of
the variables, we simply denote the elementary symmetric polynomial $S_n$ of a
corresponding zero set $B$ by $S_n(B)$.  Further, let $S_0 \coloneqq 1$ and
$S_n \coloneqq 0$ for $n<0$ and $n \ge N$.  With the help of the modified zero set and the
elementary symmetric polynomials, we can now give a criterion whether a discrete-time
signal $x$ can be uniquely recovered from its Fourier intensity $\absn{\fourier x}$ and
the modulus $\absn{x[N-1-\ell]}$ or not.

\begin{Theorem}
    \label{the:uni-mod-arb-point}%
    Let the complex-valued signal $x$ with support $\{0, \dots, N-1 \}$ be a solution of $\abs{ \fourier{x}}^{2} = \fourier{a}$
    as in (\ref{eq:char-amb}) with
    corresponding zero set $B \coloneqq \{ \beta_1, \ldots, \beta_{N-1}\}$, and
    $\ell$ be in $\{0, \ldots ,N-1\}$.  Then the signal $x$ can be uniquely recovered from
    $\absn{ \fourier x}$ and $\absn{x[N-1-\ell]}$ up to rotations if and only if
    \begin{equation}\label{sl}
        \abs{S_\ell \mleft( B \mright)}
        \ne \biggl(\,\prod_{\beta_j \in \Lambda} \abs{ \beta_j } \,\biggr)
        \cdot \abs{S_\ell \mleft( B^{(\Lambda)} \mright)} 
    \end{equation}
    holds for each non-empty subset $\Lambda \subset B$ where $\Lambda$ does not contain
    reflected zero pairs of the form
    \raisebox{0pt}[0pt][0pt]{$( \beta_j^{\,}, \overline \beta\kern0pt_j^{\,-1} )$} or zeros on the
    unit circle.
\end{Theorem}

\begin{Proof}
    By normalizing the support of the signal $x$ to $\{0, \dots, N-1 \}$, we avoid the
    trivial shift ambiguity, and the modulation $\e^{-\I \omega n_0}$ in
    Theorem~\ref{the:repr-sol-time-dom} vanishes.  Therefore, the Fourier transform of the
    original signal $x$ has the form
    \begin{equation*}
        \fourier x \mleft( \omega \mright)
        = \sum_{n=0}^{N-1} x[n] \, \e^{-\I \omega n}
        = \e^{\I \alpha} \sqrt{ \abs{ a \mleft[ N-1 \mright] } \prod_{j=1}^{N-1} \abs{
                \beta_j}^{-1}} \cdot
        \prod_{j=1}^{N-1} \left( \e^{-\I \omega} - \beta_j  \right) \!.
    \end{equation*}
    Since $\fourier x$ is an algebraic polynomial in $\e^{-\I \omega}$, Vieta's formulae
    imply
    \begin{equation*}
        \abs{ x \mleft[ N-1-\ell \mright] }
        = \sqrt{ \abs{ a \mleft[ N-1 \mright] } \prod_{j=1}^{N-1} \abs{\beta_j}^{-1}}
        \cdot \abs{ S_\ell \mleft( B \mright)}.
    \end{equation*}

    We suppose now that the signal $x$ cannot be uniquely recovered up to rotations.  Then
    there exists a second solution $\breve x$ that is not a rotation of the signal $x$.
    Since the zero sets of $x$ and $\breve x$ cannot coincide, there exists a subset
    $\Lambda$ which does not contain reflected zero pairs or zeros on the unit circle such
    that $\breve x$ corresponds to the modified zero set $B^{(\Lambda)}$.  Now, a
    comparison of the moduli $\absn{x[N-1-\ell]}$ and $\absn{\breve x [N-1-\ell]}$ yields
    \begin{equation*}
        \sqrt{ \abs{ a \mleft[ N-1 \mright] } \prod_{j=1}^{N-1} \abs{\beta_j}^{-1}}
        \cdot \abs{ S_\ell \mleft( B \mright)}
        = \sqrt{ \abs{ a \mleft[ N-1 \mright] } \prod_{j=1}^{N-1}
            \abs{\beta_j^{(\Lambda)}}^{-1}} 
        \cdot \abs{ S_\ell \mleft( B^{(\Lambda)} \mright)} \!. 
    \end{equation*}
    Simplifying this equation, we find
    \begin{equation*}
        \abs{S_\ell \mleft( B \mright)}
        = \biggl(\,\prod_{\beta_j \in \Lambda} \abs{ \beta_j } \,\biggr)
        \cdot \abs{S_\ell \mleft( B^{(\Lambda)} \mright)} \!. 
    \end{equation*}
    Hence, the unknown signal $x$ can only be recovered uniquely up to rotations if and
    only if there exists no subset $\Lambda$ fulfilling the above equation.  \qed
\end{Proof}

In the special case when the solution signal $x$ has support $\{0,\dots, N-1\}$ of odd length $N$, and we have given
the modulus of the centered value $\absn{x[\nicefrac{(N-1)}{2}]}$, then we need to pay special attention 
because $\absn{x[\nicefrac{(N-1)}{2}]}$ does not change under the
reflection and conjugation of the complete signal.  Since the reflected and conjugated
signal corresponds to the reflection of the complete zero set, see
Lemma~\ref{lem:zeros-refl-sig}, the uniqueness condition in Theorem~\ref{the:uni-mod-arb-point}
cannot be satisfied except if $B$ is invariant under the reflection at the unit circle.  To
avoid this trivial ambiguity, we assume that the second solution $\breve x$ in the proof
of Theorem~\ref{the:uni-mod-arb-point} is not equal to the reflected and conjugated signal
$x$. In other words, the subset $\Lambda$ should not be extendable to the complete set $B$
by adding reflected zero pairs or zeros on the unit circle.  Adapting the proof of
Theorem~\ref{the:uni-mod-arb-point}, we obtain the following slightly weaker
statement.

\begin{Corollary}
    \label{cor:uni-mod-arb-point}%
    Let the complex-valued signal $x$ with support $\{0,\dots, N-1\}$ of odd length $N$
    be a solution of $\abs{ \fourier{x}}^{2} = \fourier{a}$ as in (\ref{eq:char-amb}) with
    corresponding zero set $B \coloneqq \{ \beta_1, \dots, \beta_{N-1}\}$. Then the
    signal $x$ can be uniquely recovered from $\absn{ \fourier x}$ and
    $\absn{x[\nicefrac{(N-1)}{2}]}$ up to rotations and conjugate reflections if and only
    if
    \begin{equation*}
        \abs{S_{\frac{N-1}{2}} \mleft( B \mright)}
        \ne \biggl(\,\prod_{\beta_j \in \Lambda} \abs{ \beta_j } \,\biggr)
        \cdot \abs{S_{\frac{N-1}{2}} \mleft( B^{(\Lambda)} \mright)}
        \addmathskip
    \end{equation*}
    holds for each non-empty proper subset $\Lambda \subset B$ where $\Lambda$ does not
    contain reflected zero pairs or zeros on the unit circle and cannot be extended to the
    complete set $B$ by adding zeros of this kind.
\end{Corollary}

We want to show that the conditions in Theorem~\ref{the:uni-mod-arb-point} and Corollary~\ref{cor:uni-mod-arb-point}
are almost always satisfied for given measurements $\absn{ \fourier x}$ and
    $\absn{x[N-1-\ell]}$.
For this purpose, we identify the elementary
symmetric polynomial $ S_\ell \mleft( B \mright)$ with a real $(2N-2)$-variate polynomial
in the variables
\begin{equation}
    \label{eq:vec-beta}
    \Vek \beta \coloneqq
    \left( \Re \beta_1, \Im \beta_1, \dots, \Re \beta_{N-1} , \Im \beta_{N-1}
    \right)^\T
    \in \left(\R^2 \setminus \left\{0\right\}\right)^{N-1}
\end{equation}
and show that for every non-empty $\Lambda \subset B$
and each $\ell \in \{0, \dots, N-1 \}$ the exceptional zero sets are contained in the zero locus of a
non-trivial $(2N-2)$-variate algebraic polynomial.

\begin{Lemma}
    \label{lem:vari-arb-point}%
    Let $B:= \{ \beta_{1}, \ldots , \beta_{N-1} \} \subset {\C}^{N-1}$, and let $\ell$ be in $\{0, \dots, N-1\}$. 
    Then, for each non-empty subset
    $\Lambda \subset B$, the zero sets $B$ satisfying
    \begin{equation}
        \label{eq:vari-arb-point:cond}
        \abs{S_\ell \mleft( B \mright)}
        = \Bigl(\,\prod_{\beta_j \in \Lambda} \abs{ \beta_j } \,\Bigr)
        \cdot \abs{S_\ell \mleft( B^{(\Lambda)} \mright)}
    \end{equation}
    with $B^{(\Lambda)}$ in \eqref{B+} can be identified with the zero locus of a
    non-trivial polynomial in $2N-2$ variables whenever $\ell \ne
    \nicefrac{(N-1)}{2}$.
    For the case $\ell =\nicefrac{(N-1)}{2}$, the assertion holds true if $\Lambda$ is a
    proper subset of $B$.
\end{Lemma}

\begin{Proof}
    Using the definition of the elementary symmetric function \eqref{eq:ele-sys-fun}, the
    condition (\ref{eq:vari-arb-point:cond}) implies
    \begin{equation}
        \label{eq:poly-cond-alg-variety}
        \absBB{\smashoperator[r]{\sum_{1\le k_1 < \cdots < k_\ell \le N-1\hspace*{-10pt}}}
            \,\beta_{k_1} \cdots \beta_{k_\ell}}^2
        = \absBB{\prod_{\beta_j \in \Lambda} \overline \beta_j }^2 \cdot \,
        \absBB{\smashoperator[r]{\sum_{1\le k_1 < \cdots < k_\ell \le N-1\hspace*{-10pt}}}
            \,\beta_{k_1}^{(\Lambda)} \cdots \beta_{k_\ell}^{(\Lambda)}}^2 \!,
        \addmathskip
    \end{equation}
    where the empty sums for $\ell=0$ have been set to one.  With the substitution
    $\beta_j = \Re \beta_j + \I \, \Im \beta_j$, the left-hand side of
    \eqref{eq:poly-cond-alg-variety} becomes a real algebraic polynomial
    $p_{1}(\Vek \beta)$.  Since the reflection of a zero $\beta_j$ at the unit circle is
    simply given by \raisebox{0pt}[0pt][0pt]{$\overline \beta\kern0pt_j^{\,-1}$}, the
    reflected zeros in the modified zero set  $B^{(\Lambda)}$ on the right-hand side completely cancel with
    the prefactor.  Hence, the right-hand side of \eqref{eq:poly-cond-alg-variety} is also
    an algebraic polynomial $p_{2}(\Vek \beta)$.

    Thus the vectors $\Vek \beta$ satisfying \eqref{eq:poly-cond-alg-variety} form the
    zero locus of a real $(2N-2)$-variate algebraic polynomial $p_{1} - p_{2}$.  We only
    have to show that the two polynomials $p_{1}$ and $p_{2}$ on both sides of
    \eqref{eq:poly-cond-alg-variety} do not coincide.  Without loss of generality we
    assume that $\Lambda$ contains the first $J$ zeros of $B$.  Determining the real and
    imaginary parts of a summand
    \begin{equation*}
        \beta_{k_1} \cdots \beta_{k_\ell}
        = \left( \Re \beta_{k_1} + \I \Im \beta_{k_1} \right) \cdots
        \left( \Re \beta_{k_\ell} + \I \Im \beta_{k_\ell} \right)
    \end{equation*}
    on the left-hand side of \eqref{eq:poly-cond-alg-variety}, we obtain a homogeneous
    polynomial of degree $\ell$ in the real variables $\Re \beta_{k_1}$,
    $\Im \beta_{k_1}$, \dots, $\Re \beta_{k_\ell}$, $\Im \beta_{k_\ell}$.  Thus the
    polynomial $p_{1}$ is a $(2N-2)$-variate real homogeneous polynomial of degree
    $2\ell$.  By contrast, we show that $p_{2}$ is not homogeneous of degree $2\ell$ since
    it contains monomial terms of different degree.  We distinguish the following cases.

    \begin{enumerate}[(i)]
    \item For numbers $\ell$ and $J$ with $\ell + J \le N-1$, we always find increasing
        indicies $k_1< \cdots <k_\ell$ such that $k_1 > J$.  Then the product
        $\beta_{k_1}^{(\Lambda)} \cdots \beta_{k_\ell}^{(\Lambda)}$ in
        \eqref{eq:poly-cond-alg-variety} simply becomes
        $\beta_{k_1}\cdots \beta_{k_\ell}$, and hence no zeros cancel with the prefactor.
        This implies that the corresponding monomials in $p_{2}(\Vek \beta)$ are exactly
        of degree $2(\ell + J) \neq 2\ell$.
    \item If the numbers $\ell$ and $J$ fulfil $\ell + J > N-1$ and $J \le \ell$, then we
        consider the summand with the indicies $ k_1 = 1, \dots, k_\ell = \ell$.  Now, the
        first $J$ modified zeros $\beta_j^{(\Lambda)}$ cancel with the prefactor.  Since
        the real and imaginary parts of the obtained summand
        $\beta_{J+1}\cdots \beta_{\ell}$ consist of monomials of degree $\ell - J$, we
        have at least one monomial of degree $2(\ell - J) \neq \ell$ in $p_{2}$.
    \item Let finally $\ell + J > N-1$ and $J > \ell$.  We consider again the summands
        with indices $ k_1=1, \, \dots, k_\ell = \ell$ on the right-hand side of
        (\ref{eq:poly-cond-alg-variety}).  This yields the summand
        $\overline \beta_{\ell+1} \cdots \overline\beta_J$, which corresponds to monomials
        of degree $2(J-l)$ in $p_{2}(\Vek \beta)$, since all modified zeros cancel with
        the prefactor.  Thus, for $J \neq 2\ell$, these terms have a degree different from
        $2\ell$.

        For the special case $J = 2 \ell$ we consider the indices
        $ k_1 = N-\ell, \, \dots, k_\ell = N-1$ that correspond to the summand
        \begin{equation*}
            \overline \beta_1 \cdots \overline \beta_{N-\ell-1} \,
            \beta_{J+1} \cdots \beta_{N-1}
        \end{equation*}
        with $2N-2-\ell-J = 2N-2-3\ell$ different complex variables.  Thus, the
        corresponding terms of $p_{2}(\Vek \beta)$ has degree $4N-4-2\ell-2J= 4N-4-6\ell$,
        being different from $2\ell$ for $J \neq N-1$.  Hence, assuming that $J \neq N-1$,
        i.e., $\Lambda$ is a proper subset of $B$, there exists a monomial term of degree
        different from $2\ell$ also in the special case $J=2\ell$.  \qed
    \end{enumerate}
\end{Proof}

Now we can conclude the following recovery result.
\begin{Theorem}
    \label{the:almost-uni-arb-point}%
    Let $\ell$ be
    an arbitrary integer between $0$ and $N-1$.  
    Then almost every complex-valued signal $x$ with support $\{0,\dots,N-1\}$ can be
    uniquely recovered from $\absn{\fourier x}$ and $\absn{x[N-1-\ell]}$ up to rotations
    if $\ell \ne \nicefrac{(N-1)}{2}$.  In the case $\ell=\nicefrac{(N-1)}{2}$, the
    reconstruction is almost surely unique up to rotations and conjugate reflections.%
\end{Theorem}

\begin{Proof}
    From Lemma~\ref{lem:vari-arb-point} we can conclude that for all possible choices of
    $\Lambda$ and $\ell$ the exceptional zero sets which do not satisfy the uniqueness
    conditions in Theorem~\ref{the:uni-mod-arb-point} and Corollary~\ref{cor:uni-mod-arb-point} are
    contained in the union of finitely many zero loci of $(2N-2)$-variate real algebraic
    polynomials and thus form a set $E$ of Lebesgue measure zero.  It remains to show that
    this observation can be transferred to the components of the corresponding signals,
    which can be represented by
    \begin{equation}
        \label{eq:rel-zeros-comp}
        x \mleft[ N-1-n \mright] = \left(-1\right)^n C \,
        S_n \mleft( \beta_1, \dots, \beta_{N-1} \mright), \quad , n=0, \ldots , N-1,
    \end{equation}
    where $C$ is an appropriate constant.

    In order to prove this statement, we will apply the following variant of Sard's
    theorem {\cite[Theorem~3.1.]{Sch69}}.  Let $F \colon D \rightarrow \R^n$ be a
    continuously differentiable mapping where $D$ is an open set in $\R^n$. Then the image
    $F(E)$ of every measurable set $E \subset D$ is measurable, and the Lebesgue measure
    $\lambda$ of the image $F(E)$ is bounded by
    \begin{equation*}
        \lambda \mleft( F \mleft( E \mright) \mright)
        \le \int_E \abs{ \det \Mat J_F \mleft( y \mright)} \diff{y},
        \submathskip
    \end{equation*}
    where $\Mat J_F$ is the Jacobian of $F$.%

    We reconsider the equations in \eqref{eq:rel-zeros-comp} as a mapping
    $F \colon \R^{2N} \rightarrow \R^{2N}$ with
    \begin{align*}
   &F\mleft( \Re \beta_1, \Im \beta_1, \dots, \Re \beta_{N-1}, \Im \beta_{N-1},
     \Re C, \Im C \mright)
      \\[\fskip]
   &\quad = \left( \Re x\mleft[ 0 \mright], \Im x\mleft[ 0 \mright], \dots,
     \Re x\mleft[ N-1 \mright], \Im x\mleft[ N-1 \mright] \right)^\T \!.
    \end{align*}
    Due to the fact that $F$ is continuously differentiable.  We now apply Sard's theorem
    and conclude that the set $F(E \times (\R^2 \setminus \{0\}))$ of all signals which
    cannot be recovered uniquely up to rotations also have Lebesgue measure zero.
    Since the remaining signals can be uniquely reconstructed up to rotations, the
    assertion follows.  \qed%
\end{Proof}

\subsection{The moduli of the entire signal}
\label{sec:moduli-entire-signal}

Next, we will investigate the question whether every signal $x$ can be uniquely recovered
up to rotations if more then one modulus $\absn{x[n]}$ or even all moduli
$(\absn{x[n]})_{n\in\Z}$ are given.  Phase retrieval problems of this kind, where the
complete modulus of the signal is known, have been numerically studied in
\cite{GS72,SSD+06,LT08,LT09}.  Based on our findings in the last subsection, we immediately
obtain the following statement.

\begin{Corollary}
    \label{cor:almost-uni-whole-mod}%
    Almost every complex-valued signal $x$ with support $\{0,\dots,N-1\}$ 
    can be uniquely recovered from $\absn{\fourier x}$ and
    $(\absn{x[n]})_{n=0}^{N-1}$ up to rotations.%
\end{Corollary}

Unfortunately, the additional knowledge of more than one modulus of the signal in time
domain does not ensure uniqueness of the solution of the corresponding phase retrieval
problem in general, even if the complete modulus $\absn{x} \coloneqq (x[n])_{n\in\Z}$ of
the signal $x$ is given.

\begin{Theorem}
    \label{the:amb-whole-mod}%
    For every integer $N>3$, there exists a signal $x$ with support $\{0,\dots,N-1\}$ such
    that $x$ cannot be uniquely recovered from $\absn{\fourier x}$ and $\absn{x}$ up to
    rotations.%
\end{Theorem}

\begin{Proof}
    We consider the signal $x$ with support $\{0, \dots, N-1\}$ whose corresponding zero
    set is given by
    $B \coloneqq \{ \eta_1^{\,}, - \eta_1^{-1}, \I \eta_2^{\,}, \dots, \I \eta_2^{\,}\}$
    with $\eta_1, \eta_2 \in \R$ and $\eta_1, \eta_2 > 1$.  In other words, we choose
    $\beta_1 \coloneqq \eta_1^{\,}$, $\beta_2 \coloneqq - \eta_1^{-1}$, and
    $\beta_3 \coloneqq \cdots \coloneqq \beta_{N-1} = \I \eta_2^{\,}$.  Using
    Theorem~\ref{the:repr-sol-time-dom}, we can write the Fourier transform of $x$ as
    \begin{equation*}
        \fourier x(\omega)
        = \sum_{n=0}^{N-1} x[n] \e^{-\I \omega n}
        = C \left(\e^{-\I \omega} - \eta_{1}^{\,}\right)
        \left(\e^{-\I \omega} + \eta_{1}^{-1}\right)
        \left(\e^{-\I \omega} - \I \eta_2^{\,} \right)^{N-3} 
    \end{equation*}
    where $C \coloneqq \e^{\I\alpha}\sqrt{\absn{a[n-1]} \absn{\eta_2}^{N-3} }$.  Further,
    we consider the signal
    $y$ given in the frequency domain by
    \begin{equation*}
        \fourier y(\omega)
        = \sum_{n=0}^{N-1} y[n] \e^{-\I \omega n}
        = C \left(\e^{-\I \omega} - \eta_{1}^{-1} \right)
        \left(\e^{-\I \omega} + \eta_{1}^{\,} \right)
        \left(\e^{-\I \omega} - \I \eta_2^{\,} \right)^{N-3} 
    \end{equation*}
    with support $\{ 0, \ldots , N-1 \}$ and corresponding zero set
    $ \{ \eta_{1}^{-1}, -\eta_{1}^{\,}, \I \eta_2^{\,}, \ldots , \I \eta_2^{\,} \}$.  
    Obviously, we have $\absn{{\fourier x} (\omega)} = \absn{{\fourier y}(\omega)}$ since
    \begin{equation*}
        \abs{ \left(\e^{-\I \omega} - \eta_{1}^{\,} \right) \left(\e^{-\I \omega} +
                \eta_{1}^{-1} \right)}
        = \abs{ \left(\e^{-\I \omega} - \eta_{1}^{-1} \right) \left(\e^{-\I \omega} +  
                \eta_{1}^{\,} \right)}.
    \end{equation*}

    Moreover, we can show that $\abs{x[n]} = \abs{y[n]}$ for all
    $n \in \{ 0, \ldots, N-1 \}$.  Expanding the factorization of $\fourier x(\omega)$ and
    $\fourier y(\omega)$, we obtain
    \begin{align*}
      \abs{x[N-1-\ell]}
      &= \textstyle \abs{C} \abs{\binom{N-3}{\ell-2} (-1) (-\I \eta_2^{\,})^{\ell-2} 
        +  \binom{N-3}{\ell-1} (-\eta_{1}^{\,} + \eta_{1}^{-1}) (-\I \eta_2^{\,})^{\ell-1}
        +  \binom{N-3}{\ell} (-\I \eta_2^{\,})^{\ell}}
      \\[\fskip]
      &= \textstyle \abs{C} \absn{\eta_2^{\,}}^{\ell-2}
        \abs{ \binom{N-3}{\ell-2} (-1)
        +  \binom{N-3}{\ell-1} (-\eta_{1} + \eta_{1}^{-1}) (-\I \eta_2^{\,})
        +  \binom{N-3}{\ell} (-\eta_2^{2})}
    \end{align*}
    for $\ell=2, \ldots , N-3$, while
    \begin{equation*}
        \abs{y[N-1-\ell]} 
        = \textstyle \abs{C} \absn{\eta_2^{\,}}^{\ell-2}
        \abs{ \binom{N-3}{\ell-2} (-1)
            +  \binom{N-3}{\ell-1} (\eta_{1}^{\,} - \eta_{1}^{-1}) (-\I \eta_2)
            +  \binom{N-3}{\ell} (-\eta_2^{2})} \!.
    \end{equation*}
    Thus, we indeed have $ \abs{x[N-1-\ell]} = \abs{y[N-1-\ell]}$. The assertion follows
    analogously for the remaining indices $\ell=0,1,N-2, N-1$.  \qed%
\end{Proof}

\section{Using additional phase information}
\label{sec:add-phase-point}%

Now, we study the question, whether a~priori phase information about the unknown signal
$x$ can also enforce uniqueness of the solution of the discrete phase retrieval problem
$\abs{\fourier x}^2 = \fourier a$, where $\fourier a$ is the given non-negative
trigonometric polynomial of degree $N-1$.  Obviously, knowing the phase of only one
component $x[n]$ of the signal is not sufficient because of the trivial rotation
ambiguity.  Thus, we need at least the phase of two components.

Firstly, we consider the right and left end point of a signal $x$ given by
\begin{equation*}
    x\mleft[ N-1 \mright] =
    \e^{\I \alpha} \sqrt{\abs{a\mleft[N-1\mright]} \prod_{j=1}^{N-1}\abs{\beta_j}^{-1}}
    \submathskip
\end{equation*}
and
\begin{equation*}
    x \mleft[ 0 \mright] =
    \left(-1\right)^{N-1} \, \e^{\I\alpha}
    \sqrt{\abs{a\mleft[N-1\mright]} \prod_{j=1}^{N-1}\abs{\beta_j}^{-1}}
    \cdot \prod_{j=1}^{N-1} \beta_j
    \addmathskip
\end{equation*}
as characterized in Theorem~\ref{the:repr-sol-time-dom}.  The end points of all further
solutions are obtained by changing the zero set
$B \coloneqq \{ \beta_1, \dots, \beta_{N-1}\}$, where the corresponding zeros of a subset
$\Lambda \subset B$ are reflected at the unit circle.  Since the additional rotation by
$\alpha$ can be individually chosen for each ambiguity, we can assume without loss of
generality that the phase of the right end point coincides for all non-trivial solutions.
Observing that the phase of a complex number is invariant under reflection at the unit
circle, we can conclude that the phases of the left end point of all solutions are also
equal.  Therefore, knowing the phases of the two end points does not reduce the set of
non-trivial ambiguities.  Nonetheless, we will show that additional phase information for
two components of the unknown signal which are not the two end points can really enforce
almost surely uniqueness of solutions the phase retrieval problem.

\subsection{Phase of an arbitrary point and the end point}
\label{sec:phase-arbi-end}

First, we consider the special case, where we have a priori information about arg $x[N-1]$
and a further value arg $x[N-1-\ell]$ for one $\ell \in \{ 1, \ldots, N-2 \}$.  We proceed
similarly as in  section~\ref{sec:using-add-magn}.  First, we characterize the signals
that cannot be uniquely reconstructed. Then we show that the exceptional zero sets corresponding to solution ambiguities are
contained in an appropriate algebraic variety. Finally, we conclude that ambiguities can
only arise in rare special cases.  Under the assumption that the unknown signal $x$
possesses the support $\{0, \dots, N-1 \}$ and can thus be written in the form
\eqref{eq:char-amb} with $n_0 = 0$, we obtain the following uniqueness condition, where
$S_\ell$ again denotes the elementary symmetric polynomial in \eqref{eq:ele-sys-fun}.

\vspace*{-8pt}

\begin{Theorem}
    \label{the:amb-phase-end-point}%
    Let $x$ be a complex-valued signal with support $\{0,\dots,N-1\}$ and corresponding
    zero set $B \coloneqq \{\beta_1, \dots, \beta_{N-1}\}$, and let $\ell$ be an integer
    between $1$ and $N-2$.  Then the signal $x$ cannot be uniquely recovered from
    $\absn{\fourier x}$, $\arg x[N-1]$, and $\arg x[N-1-\ell]$ if and only if there exists
    a non-empty subset $\Lambda \subset B$, where $\Lambda$ does not contain reflected
    zero pairs or zeros on the unit circle, such that $B$ and $B^{(\Lambda)}$ satisfy
    \begin{equation} \label{nec}
        \Re S_\ell ( B ) \,
        \Im S_\ell ( B^{(\Lambda)} )
        - \Im S_\ell ( B ) \,
        \Re S_\ell ( B^{(\Lambda)} ) = 0
        \submathskip
    \end{equation}
    and 
    \begin{equation}
        \Re S_\ell ( B ) \, \Re S_\ell ( B^{(\Lambda)} )
        + \Im S_\ell ( B ) \, \Im S_\ell ( B^{(\Lambda)} )
        \ge 0.
        \addmathskip
    \end{equation}
\end{Theorem}

\begin{Proof}
    We assume that the phase retrieval problem to recover $x$ from $\absn{\fourier x}$,
    $\arg x[N-1-\ell]$, and $\arg x[N-1]$ possesses at least one further solution $y$.  By
    Theorem~\ref{the:repr-sol-time-dom}, there is a subset $\Lambda \subset B$ such that the
    second solution $y$ corresponds to the modified zero set $B^{(\Lambda)}$.  Since $x$
    and $y$ are different, we can assume that $\Lambda$ is non-empty and does not contain
    reflected zero pairs or zeros on the unit circle.

    We recall that the components of a signal with support $\{0,\dots,N-1\}$ are given by
    \begin{equation}
        \label{eq:repr-comp-sig}
        x \mleft[ N-1-\ell \mright] =
        \left( -1\right)^\ell  \e^{\I \alpha} \sqrt{ \abs{a\mleft[N-1\mright]}
            \prod_{j=1}^{N-1} \abs{\beta_j}^{-1}} \cdot S_\ell \mleft( B \mright)
    \end{equation}
    for $\ell \in \{ 0, \ldots, N-1 \}$ due to Vieta's formulae. For the components
    $y[N-1-\ell]$ of the ambiguity, we have an analogous representation where the
    corresponding zero set $B$ is replaced by $B^{(\Lambda)}$.  Since
    $S_{0}(B) = S_{0}(B^{(\Lambda)}) = 1$, the phases of the end points $x[N-1]$ and
    $y[N-1]$ can only coincide if $x$ and $y$ have the same rotation factor
    $\e^{\I\alpha}$.

    The second phase condition 
    $   \arg x \mleft[ N-1-\ell \mright] = \arg y \mleft[ N-1-\ell \mright]$
    is now equivalent to
    \begin{equation*}
        \arg S_\ell ( B ) = \arg S_\ell ( B^{(\Lambda)} ).
        \addmathskip
    \end{equation*}
    Thus, the  value  $S_\ell(B^{(\Lambda)})$ has to lie on the real
    ray from the origin through $S_\ell(B)$ in the complex plane, i.e.,
    \begin{equation*}
        \Re S_\ell(B) \, \Im S_\ell(B^{(\Lambda)})
        - \Im S_\ell(B) \, \Re S_\ell(B^{(\Lambda)}) = 0.
        \submathskip
    \end{equation*}
    and
    \begin{equation*}
        \Re S_\ell ( B ) \, \Re S_\ell ( B^{(\Lambda)} )
        + \Im S_\ell ( B ) \, \Im S_\ell ( B^{(\Lambda)} )
        \ge 0.
        \addmathskip
    \end{equation*}
    This assertion also holds when one or both signal values $x[N-1-\ell]$ or
    $y[N-1-\ell]$ are zero, and the corresponding phases are not uniquely defined. \qed%
\end{Proof}

We show now that the non-uniqueness condition is only rarely satisfied.
\begin{Lemma}
    \label{lem:var-phase-end-point}%
    Let $B= \{ \beta_{1}, \ldots , \beta_{N-1} \} \in {\C}^{N-1}$, and let $\ell$ be in
    $\{1, \ldots, N-2\}$.  Then for each non-empty subset $\Lambda \subset B$, the zero
    sets $B$ satisfying
    \begin{equation}
        \label{eq:var-phase-end-point:linear-equa}
        \Re S_\ell(B) \, \Im S_\ell(B^{(\Lambda)}) - \Im S_\ell(B) \, \Re S_\ell
        (B^{(\Lambda)}) = 0
    \end{equation}
    with $B^{(\Lambda)}$ in \eqref{B+} can be identified with the zero locus of a
    non-trivial polynomial in $2N-2$ variables.
\end{Lemma}

\begin{Proof}
    Again, we identify the corresponding zero set $B$ of the signal $x$ with the real
    $(2N-2)$-di\-men\-sion\-al vector $\Vek \beta$ as in \eqref{eq:vec-beta}.  With the
    substitution $\beta_j= \Re \beta_j + \I \, \Im \beta_j$, the reflected zeros at the
    unit circle are now given by
    \begin{equation*}
        \Re \overline \beta_j^{\,-1} = \tfrac{\Re \beta_j}{\left[ \Re \beta_j \right]^2 +
            \left[ \Im \beta_j \right]^2}
        \quad\text{and}\quad
        \Im \overline \beta_j^{\,-1} = \tfrac{\Im \beta_j}{\left[ \Re \beta_j \right]^2 +
            \left[ \Im \beta_j \right]^2}.
    \end{equation*}
    Thus, the  elementary symmetric polynomials
    \begin{equation*}
        S_\ell(B^{(\Lambda)}) =
        \smashoperator{\sum_{\quad 1\le k_1 < \cdots < k_n \le N-1}}
        \left( \Re \beta_{k_1}^{(\Lambda)} + \I \Im \beta_{k_1}^{(\Lambda)} \right) \cdots
        \left( \Re \beta_{k_\ell}^{(\Lambda)} + \I \Im \beta_{k_\ell}^{(\Lambda)} \right)
    \end{equation*}
    corresponding to the modified zero sets are rational $2(N-2)$-variate functions in the
    real variables $\Re \beta_{j}$ and $\Im \beta_{j}$ ($j=1, \dots, N-1$), where the
    denominator of the individual summands contains the moduli of the reflected zeros.
    Multiplying \eqref{eq:var-phase-end-point:linear-equa} with
    \begin{equation*}
        \Pi_\Lambda \coloneqq
        \prod_{\beta_j \in \Lambda} \left( \left[ \Re \beta_j \right]^2 + \left[ \Im
                \beta_j \right]^2 \right) \!,
    \end{equation*}
    we thus obtain the equivalent condition
    \begin{equation}
        \label{eq:var-phase-end-point:poly}
        \Pi_\Lambda \left( \Re S_\ell(B) \, \Im S_\ell(B^{(\Lambda)})
            -\Im S_\ell(B) \, \Re S_\ell (B^{(\Lambda)}) \right) = 0
    \end{equation}
    with an algebraic polynomial in the variables $\Re \beta_j$ and $\Im \beta_j$ on
    the left-hand side.

    In order to show that the polynomial \eqref{eq:var-phase-end-point:poly} cannot vanish
    everywhere, we use the following idea: we choose a specific monomial in the left
    summand
    \begin{equation}
        \label{eq:one-phase:left-sum}
        \Pi_\Lambda \, \Re S_\ell(B) \, \Im S_\ell(B^{(\Lambda)}),
    \end{equation}
    and show that this monomial does not occur in the right summand
    \begin{equation}
        \label{eq:one-phase:right-sum}
        \Pi_\Lambda \, \Im S_\ell (B) \, \Re S_\ell (B^{(\Lambda)}).
    \end{equation}
    Assuming that we reflect the first $J$ zeros of $B$, \ie,
    $\Lambda \coloneqq \{ \beta_1, \dots, \beta_{J} \}$, we distinguish the
    following two major cases.

    \begin{enumerate}[(i)]
    \item\label{item:var-phase-end-point:i} Firstly, we assume that $N-1>\ell \ge J \ge 1$
        and consider the specific monomial
        \begin{equation*}
            p_{1} \mleft( \Vek \beta \mright) \,
            \coloneqq \Im \beta_1  \Bigl( \prod_{k=2}^{\ell} \left[ \Re \beta_k \right]^2
            \Bigr) \,
            \Re \beta_{\ell+1}.
        \end{equation*}
        Here $p_{1}$ uniquely arises in (\ref{eq:one-phase:left-sum}) from the factor
        $\Re \beta_2 \cdots \Re \beta_{\ell+1}$ in $\Re S_{\ell}(B)$ and the factor
        ${\Im \beta_1 \, \Re \beta_2 \cdots \Re \beta_\ell}/{\Pi_{\Lambda}}$ in
        $\Im S_{\ell}(B^{(\Lambda)})$. However, $p_{1}$ is not contained in
        (\ref{eq:one-phase:right-sum}) since otherwise $\Re S_{\ell}(B^{(\Lambda)})$ has
        to contain the factor
        ${\Im \beta_1 \, \Re \beta_2 \cdots \Re \beta_\ell}/{\Pi_{\Lambda}}$ such
        that $\Pi_{\Lambda}$ cancels out. However, this is impossible because the
        nominator of the real part $\Re S_\ell(B)$ consists only of monomials with an even
        number of \bq{imaginary variables} $\Im \beta_j$.

    \item Let now $1 \le \ell \le J \le N-1$. First we assume that
        $2\ell \ge N-1$ and investigate the monomial
        \begin{equation*}
            p_{2} \mleft( \Vek \beta \mright)
            :=\Im \beta_1 \, \biggl( \prod_{k=2}^{N-1-\ell} \Re \beta_k \biggr) \,
            \biggl( \prod_{k=N-\ell}^{\ell} \left[\Re \beta_{k} \right]^2  \biggr) \,
            \biggl( \prod_{k=\ell+1}^{J} \left[\Re \beta_{k} \right]^3 \biggr) \,
            \biggl( \prod_{k=J+1}^{N-1} \Re \beta_{J+1} \biggr),
        \end{equation*}
        where the three last products can be empty. Observe that $p_{2}$ possesses the
        degree $2J$.
        Taking the definition of $\Pi_{\Lambda}$ into account, we can conclude that the
        monomial $p_{2}$ uniquely arises from $ \Re \beta_{N-\ell} \cdots \Re \beta_{N-1}$
        in $\Re S_\ell(B)$, and
        \begin{equation*}
            \frac{\Im \beta_1 \Re \beta_2 \cdots \Re \beta_\ell}{\prod\limits_{j=1}^\ell
                \left( \left[ \Re \beta_j \right]^2 + \left[ \Im \beta_j \right]^2
                \right)}
        \end{equation*}
        in $\Im S_\ell(B^{(\Lambda)})$.  Using an analogous argumentation as before, the
        monomial $p_{2}$ cannot be a term in (\ref{eq:one-phase:right-sum}) since then the
        nominator of $\Re S_\ell (B^{(\Lambda)})$ has to contain a monomial with only one
        \bq{imaginary variable} $\Im \beta_j$.

        It remains to show that the polynomial is also non-trivial for $2\ell < N-1$.
        Here we examine the monomial
        \begin{equation*}
            \Im \beta_1 \, \biggl( \prod_{k=2}^{\ell}\Re \beta_k \biggr) \,
            \biggl( \prod_{k=\ell+1}^{N-\ell-1} \left[\Re \beta_{\ell+1} \right]^2 \biggr) \,
            \biggl(\prod_{k=N-\ell}^{J} \left[\Re \beta_{k} \right]^3 \biggr) \,
            \biggl( \prod_{k=J+1}^{N-1} \Re \beta_{J+1} \biggr)
        \end{equation*}
        in the case $N-\ell \le J$ and otherwise the monomial 
        \begin{equation*}
            \Im \beta_1 \, \biggl( \prod_{k=2}^{\ell} \Re \beta_k \biggr) \,
            \biggl( \prod_{k=\ell+1}^{J} \left[\Re \beta_{\ell+1} \right]^2 \biggr) \,
            \biggl( \prod_{k=N-\ell}^{N-1} \Re \beta_{N-\ell} \biggr).
        \end{equation*}
        Again these monomials occur in (\ref{eq:one-phase:left-sum}) but not in
        (\ref{eq:one-phase:right-sum}).
        \qed
    \end{enumerate}
\end{Proof}

Considering the union of the constructed zero loci in Theorem~\ref{lem:var-phase-end-point} for
all possible subsets $\Lambda$, we can conclude that the additional phase information can
indeed enforce uniqueness of the reconstruction for almost every signal.

\begin{Theorem}
    \label{the:almost-uni-phase-end}%
    Let $\ell$ be
    an arbitrary integer between $1$ and $N-2$.  Then almost every complex-valued signal $x$ with support $\{0,\dots,N-1\}$ can be
    uniquely recovered from $\absn{\fourier x}$, $\arg x[N-1]$, and $\arg x[N-1-\ell]$.
\end{Theorem}

\begin{Proof}
    In Theorem~\ref{lem:var-phase-end-point}, we have observed that the zero sets satisfying the
    condition (\ref{nec}) in Theorem~\ref{the:amb-phase-end-point} for a specific subset
    $\Lambda$ lie in the zero locus of an algebraic polynomial.  Since this polynomial is
    non-trivial, these zero sets form a set with zero Lebesgue measure.  Using Vieta's
    formulare and Sard's theorem, we can deduce the assertion similarly as in the proof of
    Theorem~\ref{the:almost-uni-arb-point}.  \qed%
\end{Proof}

\subsection{Phase of two arbitrary points}
\label{sec:phase-two-arbi}

Let us now consider the phase retrieval problem where we have a priori information about
the phase at two inner signal points.  More precisely, we consider the recovery of an unknown
signal $x$ from its Fourier intensity $\absn{\fourier x}$ and the phases
$\arg x[N-1-\ell_{1}]$ and $\arg x[N-1-\ell_{2}]$ with
$\ell_{1}, \ell_{2} \in \{ 1, \ldots , N-2 \}$ and $\ell_{1} \neq \ell_{2}$.  Similarly as
in Theorem~\ref{the:amb-phase-end-point}, we can characterize the corresponding zero sets of all
signals that cannot be uniquely reconstructed, where we again  apply the notations of section \ref{sec:phase-retr-probl}.

\begin{Theorem}
    \label{the:amb-phase-arb-points}%
    Let $x$ be a complex-valued signal with support $\{0,\dots,N-1\}$ and corresponding
    zero set $B \coloneqq \{\beta_1, \dots, \beta_{N-1}\}$, and let $\ell_{1}$ and
    $\ell_{2}$ be different integers between $1$ and $N-2$.  Then the signal $x$ cannot be
    uniquely recovered from $\absn{\fourier x}$, $\arg x[N-1-\ell_1]$, and
    $\arg x[N-1-\ell_2]$ if and only if there exists a non-empty subset
    $\Lambda \subset B$, where $\Lambda$ does not contain reflected zero pairs or zeros on
    the unit circle, such that $B$ and $B^{(\Lambda)}$ satisfy
    \begin{align*}
      &\Re \bigl[ S_{\ell_1} (B) \bigr] \,
        \Im \Bigl[ \overline{S_{\ell_2}(B^{(\Lambda)})} \, S_{\ell_2}(B) \,
        S_{\ell_1}(B^{(\Lambda)}) \Bigr]
      \\[\fsmallskip]
      &\qquad-
        \Im \bigl[ S_{\ell_1} (B) \bigr] \,
        \Re \Bigl[ \overline{S_{\ell_2}(B^{(\Lambda)})} \, S_{\ell_2}(B) \,
        S_{\ell_1}(B^{(\Lambda)}) \Bigr]
        = 0
      \\[\subalignskip]
    \end{align*}
    and 
    \begin{align*}
      \\[\subalignskip]
   & \Re \bigl[ S_{\ell_1}(B) \bigr] \,
     \Re \Bigl[ \overline{S_{\ell_2}(B^{(\Lambda)})} \, S_{\ell_2}(B) \,
     S_{\ell_1}(B^{(\Lambda)}) \Bigr]
      \\[\fsmallskip]
   &\qquad+
     \Im \bigl[ S_{\ell_1}(B) \bigr] \,
     \Im \Bigl[ \overline{S_{\ell_2}(B^{(\Lambda)})} \, S_{\ell_2}(B) \,
     S_{\ell_1}(B^{(\Lambda)}) \Bigr]
     \ge 0.
    \end{align*}
\end{Theorem}

\begin{Proof}
    We assume that the phase retrieval problem to recover $x$ from its Fourier intensity
    and the phases $\arg x[N-1-\ell_1]$ and $\arg x[N-1-\ell_2]$ has a further solution
    $y$.  By Theorem~\ref{the:repr-sol-time-dom}, we find a subset $\Lambda \subset B$ so that
    $y$ corresponds to the modified zero set $B^{(\Lambda)}$, where $\Lambda$ does not
    contain reflected zero pairs or zeros on the unit circle.

    Recall that the components of $x$ and $y$ can be written in the form
    (\ref{eq:repr-comp-sig}), where for $y$ we replace $B$ by $B^{(\Lambda)}$ and the
    rotation factor $\e^{\I \alpha}$ by $\e^{\I \alpha_{1}}$.  Due to the trivial rotation
    ambiguity, we can always rotate the second signal $y$ such that the phases
    $\arg x[N-1-\ell_2] = \arg \left( (-1)^{\ell_{2}} \, \e^{\I \alpha} \, S_{\ell_{2}}(B)
    \right)$
    and
    $\arg y[N-1-\ell_2] = \arg \left( (-1)^{\ell_{2}} \, \e^{\I \alpha_{1}} \,
        S_{\ell_{2}}(B^{(\Lambda)}) \right)$ coincide.  In other words, we choose
    $$ \alpha_{1} = \arg \left( \e^{\I \alpha_{1}} \,
        \overline{S_{\ell_{2}}(B^{(\Lambda)})} \, S_{\ell_{2}}(B) \right). $$
    Using the representation \eqref{eq:repr-comp-sig} for 
    of $x[N-1-\ell_{1}]$ and $y[N-1-\ell_{1}]$, we can simplify the condition for the second given phase to
    \begin{equation*}
        \arg \bigl( \left(-1\right)^{\ell_1} \e^{\I \alpha} \, S_{\ell_1}(B) \bigr)
        = \arg \bigl( \left(-1 \right)^{\ell_1} \e^{\I \alpha} \, \overline{ S_{\ell_2}
            (B^{(\Lambda)})} \, S_{\ell_2}(B) \, S_{\ell_1}(B^{(\Lambda)}) \bigr)
        \addmathskip
    \end{equation*}
    and thus to
    \begin{equation*}
        \arg \mleft( S_{\ell_1}(B) \mright)
        = \arg \mleft( \overline{ S_{\ell_2} (B^{(\Lambda)})} \, S_{\ell_2}(B) \,
        S_{\ell_1}(B^{(\Lambda)}) \mright).
        \addmathskip
    \end{equation*}

    Following the procedure in the proof of Theorem~\ref{the:amb-phase-end-point}, the complex
    numbers $ S_{\ell_1}(B)$ and
    $ \overline{ S_{\ell_2} (B^{(\Lambda)})} \, S_{\ell_2}(B) \,
    S_{\ell_1}(B^{(\Lambda)})$
    have to lie on the same ray starting from the origin in the complex plane, which
    results in the linear equation and the inequality condition of the assertion.  \qed%
\end{Proof}

\begin{Remark}
    In the symmetric case, where $\arg x[\ell]$ and $\arg x [N-1-\ell]$ are given for an
    $\ell \in \{ 1, \ldots , N-2 \}$, the phase retrieval problem always yields a second
    solution $y$ of the form
    \begin{equation*}
        y \coloneqq \e^{\I \left( \arg x[\ell] + \arg x [N-1-\ell] \right)} \, \overline{x
            \mleft[ N-1-\cdot \mright]}
        \addmathskip
    \end{equation*}
    caused by rotation, shift, and conjugation and reflection of $x$. In particular, we
    obtain
    \begin{equation*}
        \arg y \mleft[ \ell \mright]
        = \arg x[\ell] + \arg x [N-1-\ell] - \arg x \mleft[ N-1-\ell \mright]
        \submathskip
    \end{equation*}
    and
    \begin{equation*}
        \arg y \mleft[ N-1-\ell \mright]
        = \arg x[\ell] + \arg x [N-1-\ell] - \arg x \mleft[ \ell \mright],
        \addmathskip
    \end{equation*}
    which implies that $y$ really is an ambiguity.  Thus, Theorem~\ref{the:amb-phase-arb-points}
    holds true in this case with $\Lambda= B$.  In order to eliminate this special case,
    we will assume that $\Lambda$ is a proper subset of $B$ whenever
    $\ell_{2} = N-1-\ell_{1}$. \qed%
\end{Remark}

\begin{Lemma}
    \label{lem:var-phase-arb-point}%
    Let $B = \{ \beta_{1}, \ldots , \beta_{N-1} \} \subset {\C}^{N-1}$, and let $\ell_1$
    and $\ell_2$ be different integers in $\{1, \dots, N-2 \}$.  Then for each non-empty
    subset $\Lambda \subset B$, the zero sets $B$ satisfying
    \begin{equation}
        \label{eq:var-phase-arb-point:cond}
        \begin{aligned}
            &\Re \bigl[ S_{\ell_1} (B) \bigr] \,
            \Im \Bigl[ \overline{S_{\ell_2}(B^{(\Lambda)})} \, S_{\ell_2}(B) \,
            S_{\ell_1}(B^{(\Lambda)}) \Bigr]
            \\[\fsmallskip]
            &\qquad-
            \Im \bigl[ S_{\ell_1} (B) \bigr] \,
            \Re \Bigl[ \overline{S_{\ell_2}(B^{(\Lambda)})} \, S_{\ell_2}(B) \,
            S_{\ell_1}(B^{(\Lambda)}) \Bigr]
            = 0
        \end{aligned}
    \end{equation}
    with $B^{(\Lambda)}$ in \eqref{B+} can be identified with the zero locus of a
    non-trivial polynomial in $2N-2$ variables whenever $\ell_1+\ell_2 \ne N-1$.  In the
    case $\ell_1 + \ell_2 = N-1$, the statement holds for every proper subset
    $\Lambda \subset B$.
\end{Lemma}

\begin{Proof}
    Similarly as in the proof of Lemma~\ref{lem:var-phase-arb-point}, we substitute
    $\beta_{j} = \Re \beta_j + \I \, \Im \beta_j$ and multiply the equation
    (\ref{eq:var-phase-arb-point:cond}) with
    \begin{equation*}
        \Pi_\Lambda \coloneqq \prod_{\beta_j \in \Lambda} \left( \left[ \Re \beta_j
            \right]^2 + \left[ \Im \beta_j \right]^2 \right)^2
    \end{equation*}
    to obtain the equivalent condition
    \begin{equation}
        \label{eq:var-phase-arb-point:poly}
        \begin{aligned}
            \Pi_\Lambda \,
            &\Bigl(\Re \bigl[ S_{\ell_1} (B) \bigr] \, \Im \Bigl[
            \overline{S_{\ell_2}(B^{(\Lambda)})} \, S_{\ell_2}(B) \,
            S_{\ell_1}(B^{(\Lambda)}) \Bigr]
            \\[\fsmallskip]
            &\qquad- \Im \bigl[ S_{\ell_1} (B) \bigr] \, \Re \Bigl[
            \overline{S_{\ell_2}(B^{(\Lambda)})} \, S_{\ell_2}(B) \,
            S_{\ell_1}(B^{(\Lambda)}) \Bigr] \Bigr) = 0,
        \end{aligned}
    \end{equation}
    whose left-hand side can be understood as an $2(N-1)$-variate algebraic polynomial in
    the real variables $\Re \beta_j$ and $\Im \beta_j$.  We show that this polynomial is
    non-trivial by finding in any case a monomial that is contained in
    \begin{equation}
        \label{eq:var-phase-arb-point:right-poly}
        \Pi_\Lambda \, \Im \bigl[ S_{\ell_1} (B) \bigr] \, \Re \Bigl[
        \overline{S_{\ell_2}(B^{(\Lambda)})} \, S_{\ell_2}(B) \, S_{\ell_1} (B^{(\Lambda)}) \Bigr]
    \end{equation}
    but not in
    \begin{equation}
        \label{eq:var-phase-arb-point:left-poly}
        \Pi_\Lambda \, \Re \bigl[ S_{\ell_1} (B) \bigr] \, \Im \Bigl[
        \overline{S_{\ell_2}(B^{(\Lambda)})} \, S_{\ell_2}(B) \, S_{\ell_1}
        (B^{(\Lambda)}) \Bigr].
        \addmathskip
    \end{equation}
    Since the polynomial in \eqref{eq:var-phase-arb-point:poly} depends on $\ell_1$,
    $\ell_2$, and $J$, this leads to a cumbersome case study.  We discuss only one case in
    detail and state for all other cases the monomials that can be shown to arise in
    (\ref{eq:var-phase-arb-point:right-poly}) but not in
    (\ref{eq:var-phase-arb-point:left-poly}).

    \begin{enumerate}[(i)]
    \item\label{item:var-phase-arb-point:i} Let $\ell_1 > \ell_2 \ge J$ and consider the monomial 
        \begin{equation}
            \label{eq:var-phase-arb-point:mono:1}
            \textstyle  { [ \Re \beta_1 ]^2 \,
                \Bigl( \prod\limits_{k=2}^{\ell_{2}} [ \Re \beta_k ]^4 \Bigr) \,
                [ \Re \beta_{\ell_2+1} ]^3 \,
                \Bigl( \prod\limits_{k=\ell_{2}+2}^{\ell_{1}} [ \Re \beta_{k} ]^2 \Bigr)
                \, \Im \beta_{\ell_1+1},}
        \end{equation}
        where the product over the squared variables can be empty.  
        The degree of this monomial is $2 ( \ell_1 + \ell_2 )$.  Consequently,
        the monomial in (\ref{eq:var-phase-arb-point:mono:1}) is obtained from the term
        $   \Re \beta_2 \cdots \Re \beta_{\ell_1} \Im \beta_{\ell_1+1}$
        in $\Im [ S_{\ell_1}(B) ]$ and
        \begin{equation*}
            \Pi_\Lambda^{-1} \,  [ \Re \beta_1 ]^2 \, [ \Re \beta_2 ]^3 \cdots [ \Re
                \beta_k ]^3 \, [ \Re \beta_{k+1} ]^2  \, \Re \beta_{k+2} \cdots
                \Re \beta_{\ell}
        \end{equation*}
        in
        $\Re [\overline{S_{\ell_2}(B^{(\Lambda)})} \, S_{\ell_2}(B) \, S_{\ell_1}
        (B^{(\Lambda)})]$.
        Since this factorization is unique, the considered monomial
        \eqref{eq:var-phase-arb-point:mono:1} cannot vanish within
        \eqref{eq:var-phase-arb-point:right-poly}.
        Recalling that  monomials in the
        real part of $S_{\ell_1}(B)$ contain an even number of \bq{imaginary variables}
        $\Im \beta_j$, we can simply conclude, that the monomial (\ref{eq:var-phase-arb-point:mono:1})
        does not occur in (\ref{eq:var-phase-arb-point:left-poly}).

    \item\label{item:var-phase-arb-point:ii} $\ell_1 \ge J > \ell_2$ with  $2\ell_2 \le J$:
        \begin{equation*}
            \textstyle      [ \Re \beta_1 ]^2 \,
            \Bigl( \prod\limits_{k=2}^{\ell_{2}}   [ \Re \beta_k ]^3 \Bigr) \,
            \Bigl( \prod\limits_{k=\ell_{2}+1}^{J-\ell_{2}}  [ \Re \beta_{k} ]^4 \Bigr) \,
            \Bigl( \prod\limits_{k=J-\ell_{2}+1}^{J}   [ \Re \beta_{k} ]^5 \Bigr) \,
            \Bigl( \prod\limits_{k=J+1}^{\ell_{1}}  [ \Re \beta_{k} ]^2 \Bigr) 
            \, \Im \beta_{\ell_1+1}
        \end{equation*}

    \item\label{item:var-phase-arb-point:iii} $\ell_1 \ge J > \ell_2$ with $2\ell_2 > J$:
        \begin{equation*}
            \textstyle    [ \Re \beta_1 ]^2 \,
            \Bigl( \prod\limits_{k=2}^{J-\ell_{2}}  [ \Re \beta_k ]^3 \Bigr) \,
            \Bigl( \prod\limits_{k=J-\ell_{2}+1}^{\ell_{2}} [ \Re \beta_{k} ]^4 \Bigr) \,
            \Bigl( \prod\limits_{k=\ell_{2}+1}^{J} [ \Re \beta_{k} ]^5 \Bigr) \,
            \Bigl( \prod\limits_{k=J+1}^{\ell_{1}} [ \Re \beta_{k} ]^2 \Bigr)
            \, \Im \beta_{\ell_1+1}
        \end{equation*}

    \item\label{item:var-phase-arb-point:iv} $\ell_2 < \ell_1 < J$
        with $2 \ell_2 < 2 \ell_1 \le J$:
        \begin{equation*}
            \begin{aligned}
                \textstyle  \Bigl( \prod\limits_{k=1}^{\ell_{2}} [ \Re \beta_k ]^2  \Bigr) \,
                \Bigl( \prod\limits_{k=\ell_{2}+1}^{\ell_{1}} [ \Re \beta_{k} ]^3 \Bigr) \,
                \Bigl( \prod\limits_{k=\ell_{1}+1}^{J-\ell_{1}} [ \Re \beta_{k} ]^4 \Bigr) \,
                [\Im \beta_{J-\ell_1+1} ]^5 &  \\[\fsmallskip]
                \cdot \,
                \textstyle \Bigl( \prod\limits_{k=J-\ell_{1}+2}^{J-\ell_{2}} [\Re \beta_{k} ]^5 \Bigr) \,
                \Bigl( \prod\limits_{k=J-\ell_{2}+1}^{J}[ \Re \beta_{k} ]^6 \Bigr) &
            \end{aligned}
        \end{equation*}

    \item\label{item:var-phase-arb-point:v}  $\ell_2 < \ell_1 < J$
        with $2 \ell_2 \le J < 2\ell_1$ and $\ell_1 + \ell_2 < J$:
        \begin{equation*}
            \begin{aligned}
                \textstyle     \Bigl( \prod\limits_{k=1}^{\ell_{2}}   [ \Re \beta_k ]^2 \Bigr) \,
                \Bigl( \prod\limits_{k=\ell_{2}+1}^{J-\ell_{1}}  [ \Re \beta_{k} ]^3 \Bigr) \,
                \Bigl( \prod\limits_{k=J-\ell_{1}+1}^{\ell_{1}}    [ \Re \beta_{k} ]^4 \Bigr) &
                \\[\fsmallskip]
                \cdot \,
                \textstyle      [ \Im \beta_{\ell_1+1} ]^5 \,
                \Bigl( \prod\limits_{k=\ell_{1}+2}^{J-\ell_{2}} [ \Re \beta_{k} ]^5 \Bigr) \,
                \Bigl( \prod\limits_{k=J-\ell_{2}+1}^{J}[ \Re \beta_{k} ]^6 \Bigr) &
            \end{aligned}
        \end{equation*}

    \item\label{item:var-phase-arb-point:vi} $\ell_2 < \ell_1 < J$ with
        $2\ell_2 \le J<2 \ell_1$ and $\ell_1 + \ell_2 = J < N-1$:
        \begin{equation*}
            \textstyle \Bigl( \prod\limits_{k=1}^{\ell_{2}}  [ \Re \beta_k ]^2 \Bigr) \,
            [ \Re \beta_{\ell_2+1} ]^3 \,
            \Bigl( \prod\limits_{k=\ell_{2}+2}^{\ell_{1}} [ \Re \beta_{k} ]^4 \Bigr) \,
            \Bigl( \prod\limits_{k=\ell_{1}+1}^{J}    [ \Re \beta_{k} ]^6 \Bigr)
            \, \Im \beta_{J+1}
        \end{equation*}

    \item\label{item:var-phase-arb-point:vii} $\ell_2 < \ell_1 < J$ with $2\ell_2 \le J
        <2\ell_1$ and $\ell_1 + \ell_2 > J$:
        \begin{equation*}
            \begin{aligned}
                & \textstyle \Bigl( \prod\limits_{k=1}^{J-\ell_{1}} [ \Re \beta_k ]^2 \Bigr) \,
                [ \Im \beta_{J-\ell_1+1} ]^3 \,
                \Bigl( \prod\limits_{k=J-\ell_{1}+2}^{\ell_{2}} [ \Re \beta_{k} ]^3 \Bigr) \,
                \Bigl( \prod\limits_{k=\ell_{2} +1}^{J-\ell_{2}} [ \Re \beta_{k} ]^4 \Bigr) \,
                \\[\fsmallskip]
                & \qquad \cdot 
                \textstyle   \Bigl( \prod\limits_{k=J-\ell_{2}+1}^{\ell_{1}} [ \Re \beta_{k} ]^5 \Bigr) \,
                \Bigl( \prod\limits_{k=\ell_{1}+1}^{J} [ \Re \beta_{k} ]^6 \Bigr)
            \end{aligned}
        \end{equation*}

    \item\label{item:var-phase-arb-point:viii} $\ell_2 < \ell_1 < J$ with $2\ell_1 > 2\ell_2 > J$:
        \begin{equation*}
            \begin{aligned}
                & \textstyle 
                \Bigl( \prod\limits_{k=1}^{J-\ell_{1}}[ \Re \beta_k ]^2 \Bigr) \,
                [ \Im \beta_{J-\ell_1+1} ]^3 \,
                \Bigl( \prod\limits_{k=J-\ell_{1}+2}^{J-\ell_{2}} [ \Re \beta_{k} ]^3 \Bigr) \,
                \Bigl( \prod\limits_{k=J-\ell_{2}+1}^{\ell_{2}} [ \Re \beta_{k} ]^4 \Bigr) \,
                \\[\fsmallskip]
                & \qquad \cdot
                \textstyle \Bigl( \prod\limits_{k=\ell_{2}+1}^{\ell_{1}} [ \Re \beta_{k} ]^5 \Bigr) \,
                \Bigl( \prod\limits_{k=\ell_{1}+1}^{J} [ \Re \beta_{k} ]^6 \Bigr)
            \end{aligned}
        \end{equation*}
    \end{enumerate}

    Summarizing, one can show that the polynomial on the left-hand side of
    \eqref{eq:var-phase-arb-point:poly} is non-trivial for all possible combinations of
    $\ell_1$, $\ell_2$, and $J$. 
    \qed%
\end{Proof}

As before, Lemma~\ref{lem:var-phase-arb-point} implies now the almost sure uniqueness of the solutions
of the discrete phase retrieval problem.

\begin{Theorem}
    \label{the:almost-uni-phase-arb}%
    Let $\ell_1$ and $\ell_2$ be different integers between $1$ and
    $N-2$.  Then almost every signal complex-valued $x$ with support 
    $\{0,\dots,N-1\}$ can be uniquely recovered from $\absn{\fourier x}$,
    $\arg x[N-1-\ell_1]$, and $\arg x[N-1-\ell_2]$ whenever $\ell_1 + \ell_2 \ne N-1$.  In the case 
    $\ell_1+\ell_2=N-1$, the recovery is only unique up to conjugate reflections.%
\end{Theorem}

\begin{Proof}
    The assertion follows in an analogous way as Theorem~\ref{the:almost-uni-phase-end}. Again,
    the corresponding zero sets satisfying the non-uniqueness conditions in
    Theorem~\ref{the:amb-phase-arb-points} 
    for a specific $\ell_1$, $\ell_2$, and $\Lambda$ lie in the zero locus of an algebraic
    polynomial by Lemma~\ref{lem:var-phase-end-point}.  Due to the fact that there exist only
    finitely many different subsets $\Lambda$, the exceptional zero sets of signals
    without a unique reconstruction form a set with zero Lebesgue measure.  With Vieta's
    formulae and Sard's theorem, the assertion follows.  \qed%
\end{Proof}

\subsection{The phase of the entire signal}
\label{sec:phase-entire-signal}

Finally, we will consider the question whether every signal $x$ can be uniquely recovered
from $\abs{\fourier{x}}$ and the complete phase information
$\arg x \coloneqq (\arg x[n])_{n =0}^{N-1}$ in the time domain. Based on
Theorem~\ref{the:almost-uni-phase-end} and Theorem~\ref{the:almost-uni-phase-arb} we obviously have
the following statement.

\begin{Corollary}
    \label{cor:almost-uni-whole-phase}%
    Almost every signal complex-valued $x$ with support $\{0,\dots,N-1\}$ can be uniquely 
    recovered from $\absn{\fourier x}$ and $\arg x$.
\end{Corollary}

Unfortunately, the complete phase information of a signal fails to enforce the 
uniqueness of the phase retrieval problem for every signal.

\begin{Theorem}
    \label{the:amb-whole-phase}%
    For every integer $N > 2$, there exists a signal $x$ with support $\{0,\dots,N-1\}$
    such that $x$ cannot be uniquely recovered from $\absn{\fourier x}$ and $\arg x$.%
\end{Theorem}

\begin{Proof}
    Let $x$ be a signal with support $\{0, \dots, N-1 \}$
    of the form 
    \begin{equation*}
        \fourier x \mleft( \omega \mright) = \sum_{n=0}^{N-1} x[n] \e^{-\I \omega n} 
        = \sqrt{\abs{a \mleft[ N-1 \mright]} \prod_{j=1}^{N-1} \abs{\beta_j}^{-1}} \cdot
        \prod_{j=1}^{N-1} \left( \e^{-\I \omega} - \beta_j \right) \!.
    \end{equation*}
    where all zeros $\beta_{j}$ are real and negative.  Since the
    linear factors $(\e^{-\I\omega} - \beta_j)$ only have positive coefficients, the
    components $x[n]$ of the constructed signal $x$ are real and non-negative, i.e.
    $\arg x[n] =0$ for all $n$.  This observation remains valid for all arising
    ambiguities in Theorem~\ref{the:repr-sol-time-dom}.  Hence, if the corresponding zero set
    contains at least two zeros unequal to $-1$, we find a further ambiguity with phase
    zero in the time domain.  \qed%
\end{Proof}

\section*{Acknowledgements}

We gratefully acknowledge the funding of this work by the DFG in the frame work of the
SFB~755 \bq{Nanoscale photonic imaging} and of the GRK~2088 \bq{Discovering structure in
    complex data: Statistics meets Optimization and Inverse Problems.}

\bibliographystyle{alphadinUK}
{\footnotesize \bibliography{LITERATURE}}

\end{document}